\documentclass{article}
\usepackage[top=10truemm,left=20truemm,right=20truemm,bottom=10truemm,includefoot]{geometry}
\usepackage{graphicx}
\usepackage{amsmath,amssymb}
\usepackage{cite}
\usepackage[hidelinks]{hyperref}
\usepackage{authblk}

\title{Revisiting the Brunner-Munzel test from the viewpoint of\\local linear approximation}
\author[1]{Makito Oku}
\affil[1]{Research Center for Pre-Disease Science, University of Toyama, Japan}
\date{September 25, 2026}

\begin{document}
\maketitle

\section*{Abstract}
The Brunner-Munzel (BM) test is a nonparametric test for two independent samples that evaluates whether observations from one group tend to be greater than observations from another group, or vice versa. The BM test has a broader scope of application than the Mann-Whitney $U$ test because it does not assume equal variances between the two groups. However, the meaning of the BM test statistic is difficult to understand intuitively, which may be one of the factors hindering the widespread use of the BM test. To alleviate this problem, in this paper, I introduce an alternative interpretation of the BM test statistic from the viewpoint of local linear approximation. It is shown that the variance estimator for the sample stochastic superiority used in the BM test can be derived using local linear approximation, in which the influence of each observation on the sample stochastic superiority is assumed to be additive. This simple interpretation will help practitioners decide to use the BM test without hesitation.

\section{Introduction}

Nonparametric tests, such as the Wilcoxon signed-rank test \cite{wilcoxon1945}, the Mann-Whitney $U$ test \cite{mann1947}, and the Kruskal-Wallis test \cite{kruskal1952}, are widely used when observations from each group do not follow a normal distribution. The $U$ test is commonly used to compare two independent samples. It evaluates whether observations from one group tend to be greater than observations from another group, or vice versa. However, the $U$ test assumes that the two distributions have the same shape, and its output becomes unreliable especially when both the variances and sample sizes are unequal between the two groups \cite{pratt1964, murphy1976, penfield1994}.

The Brunner-Munzel (BM) test \cite{brunner2000} is another nonparametric test for two independent samples. It can be seen as a generalization of the $U$ test because it does not assume equal variances between the two groups. It also allows for differences in skewness and higher moments. The relationship between the $U$ test and the BM test is similar to that of the Student's $t$-test \cite{student1908} and the Welch's $t$-test \cite{welch1947}. Since the BM test has a broader scope of application than the $U$ test, some researchers recommend to use the BM test instead of the $U$ test when a nonparametric test is really needed---except for the case where the sample size of each group is large enough to justify the application of the Welch's $t$-test based on the central limit theorem---and the variances of the two groups differ largely \cite{natori2014, noguchi2021, karch2021, karch2023}.

Unfortunately, the BM test is not as common as the $U$ test. A possible reason is that the meaning of the BM test statistic is difficult to understand intuitively. The test statistic is defined using in-group ranks and pooled ranks. It is not immediately apparent how these ranks are related to the \textit{stochastic superiority}, that is, the probability of a randomly selected observation from one group is higher than a randomly selected observation from another group. The complex definition of the BM test statistic may be one of the factors hindering the widespread use of the BM test.

To alleviate this problem, in this paper, I introduce an alternative interpretation of the BM test statistic from the viewpoint of local linear approximation. It is shown that the variance estimator for the sample stochastic superiority used in the BM test can be derived using local linear approximation, in which the influence of each observation on the sample stochastic superiority is assumed to be additive. This corresponds to varying each variable one by one while keeping the other variables fixed. This simple interpretation will help practitioners decide to use the BM test without hesitation.

It should be noted that recently another nonparametric test called the $C^2$ test (or C-square test) \cite{schuurhuis2025} has been proposed. It outperforms the BM test in broad conditions. This paper does not intend to hinder the adoption of the $C^2$ test. Since it may take some time for the $C^2$ test to be implemented in major statistical tools, it would be beneficial to encourage the use of the BM test for the time being.

The reminder of this paper is organized as follows. Section 2 outlines the BM test. Section 3 outlines the $C^2$ test. Section 4 outlines the local linear approximation method. Section 5 introduces the alternative interpretation of the BM test. Section 6 discusses the relation to jackknife resampling. Section 7 presents the conclusions.

\section{The Brunner-Munzel test}

The BM test \cite{brunner2000} is a nonparametric test for comparing two independent samples $x_1,\ldots,x_n\sim F_x$ and $y_1,\ldots,y_m\sim F_y$, where $F_x$ and $F_y$ denote distributions that are not assumed to be normal and may have different variances. $n$ and $m$ denote the sample sizes of the two groups. The stochastic superiority $\theta$ is defined as follows:
\begin{equation}
\theta=P(X<Y)+\frac{1}{2}P(X=Y),
\end{equation}
where $X\sim F_x$ and $Y\sim F_y$ are two random variables that follow their respective distributions. Even if there are ties, $\theta=1/2$ means $P(X<Y)=P(X>Y)$. In the case of a two-tailed test, the null hypothesis corresponds to $\theta=1/2$, and the alternative hypothesis corresponds to $\theta\neq 1/2$. To avoid confusion with the $p$-value, the conventional notation $p$ for the stochastic superiority is not used in this paper.

\begin{figure}[t]
\centering
\includegraphics[width=0.8\linewidth]{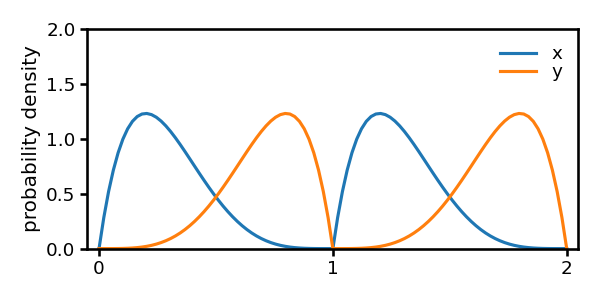}
\caption{Example case where two distributions have the same median but the stochastic superiority is not $1/2$.}
\label{fig1}
\end{figure}

The BM test is not intended to compare the medians of the two groups. Figure~\ref{fig1} shows an example case where the two distributions have the same median but $\theta\neq 1/2$. A beta distribution $B(2,5)$ defined on the interval $[0,1]$ is repeated on the interval $[1,2]$ for $x$-distribution. Another beta distribution $B(5,2)$ defined on the interval $[0,1]$ is repeated on the interval $[1,2]$ for $y$-distribution. Their densities are halved so that the integral equals to 1. By construction, the median is 1 for both distributions. There are four cases: (1) If $X<1$ and $Y>1$, $P(X<Y)=1$ holds; (2) if $X>1$ and $Y<1$, $P(X<Y)=0$ holds; (3) if $X<1$ and $Y<1$, $P(X<Y)>1/2$ holds; and (4) if $X>1$ and $Y>1$, $P(X<Y)>1/2$ holds. Since the four cases occur with equal probability, the overall tendency is $\theta>1/2$. In general, the statistical significance detected by the BM test does not imply a difference in medians.

Generally, $\theta$ can be expressed using the Riemann-Stieltjes integral as follows \cite{brunner2018}:
\begin{equation}
\theta=\int F_x\mathrm{d}F_y.
\end{equation}

\begin{figure}[t]
\centering
\includegraphics[width=0.6\linewidth]{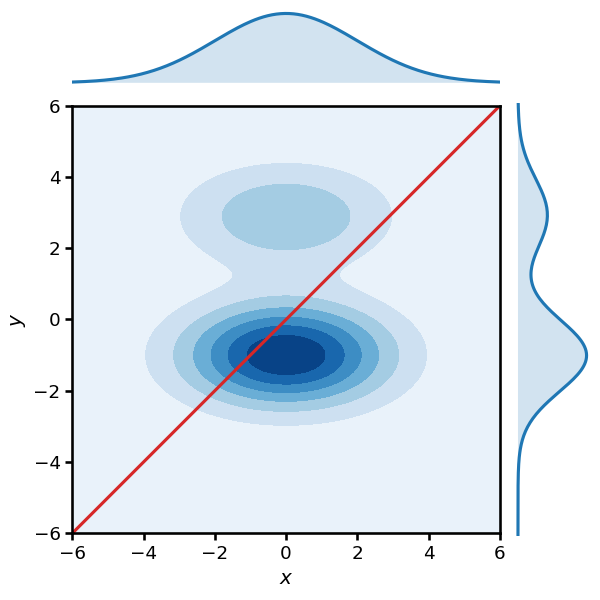}
\caption{Example of the joint probability distribution of two independent random variables. The red line represents the line $y=x$.}
\label{fig2}
\end{figure}

This expression is useful for theoretical analysis because it can be applied to both discrete and continuous distributions. If both distributions are continuous, $\theta$ can be rewritten using the Riemann integral as follows:
\begin{align}
\theta&=\int_{-\infty}^\infty F_x(y)f_Y(y)\mathrm{d}y,\\
&=\int_{-\infty}^\infty\left(\int_{-\infty}^yf_X(x)\mathrm{d}x\right)f_Y(y)\mathrm{d}y,\\
&=\iint_{x<y}f_X(x)f_Y(y)\mathrm{d}x\mathrm{d}y,\label{eq_theta_int}
\end{align}
where $f_X(x)$ and $f_Y(y)$ denote respective probability density functions. When investigating the behavior of the BM test under the null hypothesis using numerical simulations, it is necessary to make the integral value as close to 0.5 as possible, rather than matching the medians of the two distributions. Figure~\ref{fig2} shows an example of the joint probability distribution of two independent random variables: $X\sim\mathcal{N}(0,4)$ and $Y\sim 0.7\times\mathcal{N}(-1,1)+0.3\times\mathcal{N}(2.9,1)$. The integral (\ref{eq_theta_int}) corresponds to the triangular region above the line $y=x$. In this case, $\theta$ is very close to 0.5. Notice that the joint probability density function $f_{X,Y}(x,y)$ satisfies $f_{X,Y}(x,y)=f_X(x)f_Y(y)$ because $X$ and $Y$ are independent. Therefore, selecting a point $(X,Y)$ on the two dimensional plane according to its density is equivalent to selecting $X$ and $Y$ separately according to their respective densities. In short, the BM test evaluates whether the $y=x$ line bisects the joint distribution, if $X$ and $Y$ are continuous variables. That is the reason why the BM test is robust against outliers.

The sample stochastic superiority $\hat\theta$ is written by using the Heaviside function $H(x)$ as follows:
\begin{align}
\hat\theta&=\frac{1}{nm}\sum_{i=1}^n\sum_{j=1}^m H(y_j-x_i),\label{eq_theta_hat}\\
H(v)&=\left\{\begin{array}{cl}1&\text{if}\quad v>0,\\1/2&\text{if}\quad v=0,\\0&\text{if}\quad v<0.\end{array}\right.
\end{align}

\begin{figure}[t]
\centering
\includegraphics[width=0.6\linewidth]{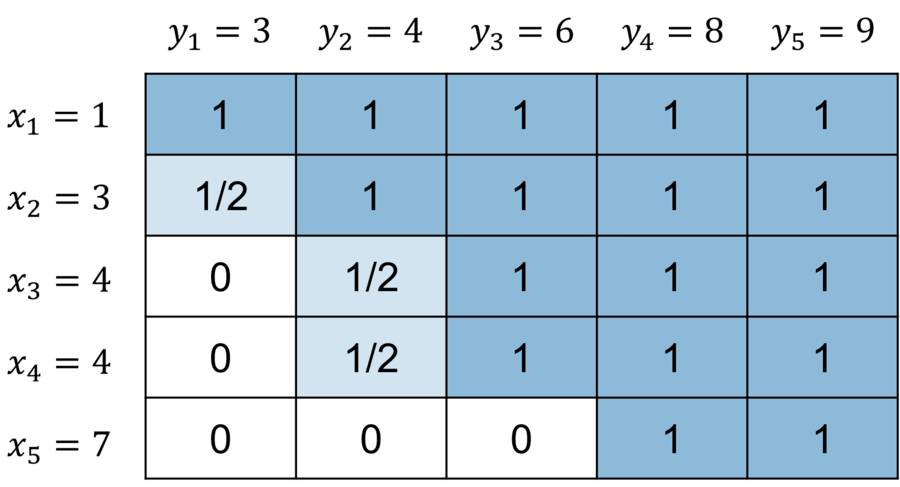}
\caption{Example data for calculating the sample stochastic superiority.}
\label{fig3}
\end{figure}

\begin{figure}[t]
\centering
\includegraphics[width=0.5\linewidth]{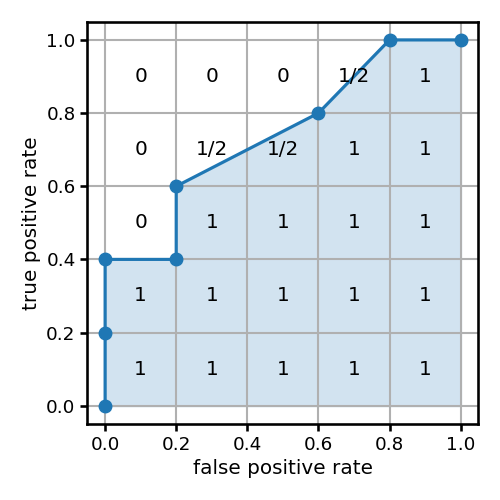}
\caption{Receiver operating characteristic curve corresponding to Fig.~\ref{fig3}.}
\label{fig4}
\end{figure}

Figure~\ref{fig3} shows an example data for calculating $\hat\theta$. Observations are assumed to be sorted within each group for simplicity. All pairs of $(x_i,y_j)$ are considered, and $\hat\theta$ represents the balance between zeros and ones. It is also known that $\hat\theta$ corresponds to the area under the receiver operating characteristic (ROC) curve \cite{delong1988}. It can be seen that rotating this table 90 degrees clockwise yields an ROC curve (Fig.~\ref{fig4}); the bottom right corner of the table becomes (0,0) point of the ROC curve, and the top left corner becomes (1,1) point. A simple interpretation is that $x_i$ and $y_i$ correspond to negative and positive cases, respectively, and their values are used for binary classification.


\begin{figure}[t]
\centering
\includegraphics[width=0.8\linewidth]{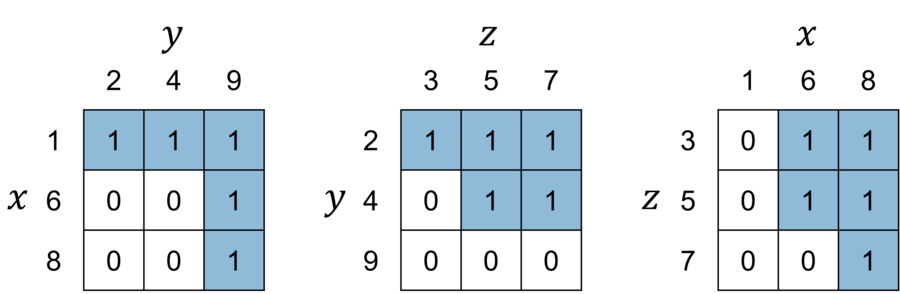}
\caption{Example of non-transitive case.}
\label{fig5}
\end{figure}

Notice that the relation based on $\hat\theta$ is not transitive \cite{thangavelu2007}. For example, suppose that there are three groups, each consisting three observations: $x=(1,6,8)$, $y=(2,4,9)$, and $z=(3,5,7)$. Let $x\prec y$ denote the condition $\hat\theta>1/2$ in the comparison between $x$ and $y$. Then, we see that $x\prec y$, $y\prec z$, and $z\prec x$, as shown in Fig.~\ref{fig5}.

Although (\ref{eq_theta_hat}) is easy to understand intuitively, the BM test statistic is defined using in-group ranks and pooled ranks. Let $r_i^x$ denote the in-group rank of $x_i$ within $x_1,\ldots,x_n$. Similarly, let $r_i^y$ denote the in-group rank of $y_i$ within $y_1,\ldots,y_m$. In addition, let $R_i^x$ and $R_i^y$ denote pooled ranks of $x_i$ and $y_i$, respectively, within the pooled observations $x_1,\ldots,x_n,y_1,\ldots,y_m$. If there are ties, midranks are used. These ranks satisfy the following relations:
\begin{equation}
\sum_{i=1}^n r_i^x=\frac{n(n+1)}{2},\quad \sum_{i=1}^m r_i^y=\frac{m(m+1)}{2},\quad \sum_{i=1}^n R_i^x+\sum_{i=1}^m R_i^y=\frac{N(N+1)}{2},
\end{equation}
where $N=n+m$ is the total number of observations.

Let $R_x$ and $R_y$ denote the sum of $R_i^x$ and $R_i^y$, respectively. Similarly, let $\bar R_x$ and $\bar R_y$ denote the mean of $R_i^x$ and $R_i^y$, respectively. 
\begin{equation}
R_x=\sum_{i=1}^n R_i^x,\quad R_y=\sum_{i=1}^m R_i^y,\quad \bar R_x=\frac{R_x}{n},\quad \bar R_y=\frac{R_y}{m}.
\end{equation}

In the BM test, $\hat\theta$ is rewritten as follows:
\begin{equation}
\hat\theta=\frac{\bar R_y-\bar R_x}{N}+\frac{1}{2}.\label{eq_theta_hat2}
\end{equation}

It is not immediately apparent how this expression is related to the sample stochastic superiority defined in (\ref{eq_theta_hat}).
To check this, we need to use the following relation:
\begin{equation}
\sum_{i=1}^n H(y_j-x_i)=R_j^y-r_j^y,\quad j=1,\ldots,m.\label{eq_rank_diff}
\end{equation}

\begin{figure}
\centering
\includegraphics[width=0.8\linewidth]{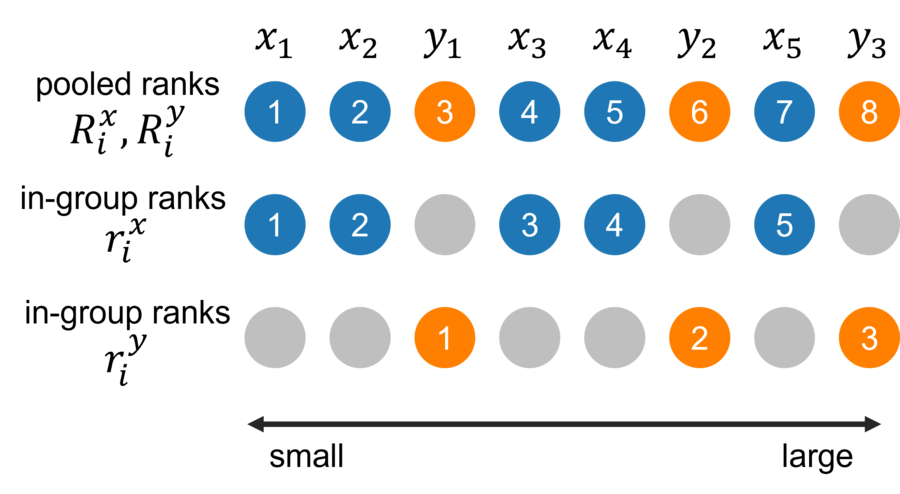}
\caption{Example case of pooled ranks and in-group ranks. Blue circles represent observations of $x_i$. Orange circles represent observations of $y_i$. All observations are sorted according to the pooled ranks. Numbers of the top row represent the pooled ranks of $x_i$ or $y_i$. Numbers of the middle row represent the in-group ranks of $x_i$. Numbers of the bottom row represent the in-group ranks of $y_i$.}
\label{fig6}
\end{figure}

This means that, for each $y_j$, the difference between the pooled rank $R_j^y$ and the in-group rank $r_j^y$ is equal to the number of observations in the other group that are smaller than $y_j$, if there is no tie. In fact, this relation holds true even when there are ties. Figure~\ref{fig6} shows an example case of pooled ranks and in-group ranks. As for $y_1$, $R_1^y-r_1^y=3-1=2$, and only $x_1$ and $x_2$ are smaller than $y_1$. As for $y_2$, $R_2^y-r_2^y=6-2=4$, and actually $x_1,\ldots,x_4$ are smaller than $y_2$. As for $y_3$,  $R_3^y-r_3^y=8-3=5$, and all 5 observations of $x_i$ are smaller than $y_3$. The reason can be understood by comparing the cases when both ranks increase against the cases when only the pooled rank increases.

The $U$ test assumes that $F_x=F_y$ holds. In this case, $R_i^x$ and $R_i^y$ can take all values with equal probability. Since $U=nm\hat\theta=\sum_{j=1}^m(R_j^y-r_j^y)=R_y-m(m+1)/2$, we get $\mathrm{Var}(U)=\mathrm{Var}(R_y)$. In the $U$ test, $\mathrm{Var}(R_y)$ is estimated under the assumption that $m$ points are randomly selected from the finite set $\{1,\ldots,N\}$ without replacement.

On the other hand, if $F_x\neq F_y$, the distributions of $R_i^x$ and $R_i^y$ become non-uniform. For example, if $F_x$ is a uniform distribution on the interval $[-1,1]$, and if $F_y$ is a uniform distribution on the interval $[-100,100]$, $R_i^x$ should be concentrated near the center. Therefore, the $U$ test should be avoided when the shapes of the two distributions differ largely because the variance estimator for $U$ becomes inaccurate.

By using (\ref{eq_rank_diff}), it can be shown that (\ref{eq_theta_hat}) and (\ref{eq_theta_hat2}) are equivalent as follows.
\begin{align}
\hat\theta&=\frac{1}{nm}\sum_{j=1}^m\sum_{i=1}^nH(y_j-x_i),\\
&=\frac{1}{nm}\sum_{j=1}^m(R_j^y-r_j^y),\\
&=\frac{1}{nm}\left(R_y-\frac{m(m+1)}{2}\right),\\
&=\frac{R_y}{nm}-\frac{m+1}{2n}-\frac{n}{2n}+\frac{1}{2},\\
&=\frac{R_y}{nm}-\frac{N+1}{2n}+\frac{1}{2},\\
&=\frac{R_y}{nm}-\frac{1}{nN}\cdot\frac{N(N+1)}{2}+\frac{1}{2},\\
&=\frac{R_y}{nm}-\frac{R_x+R_y}{nN}+\frac{1}{2},\\
&=\left(\frac{1}{nm}-\frac{1}{nN}\right)R_y-\frac{R_x}{nN}+\frac{1}{2},\\
&=\frac{R_y}{mN}-\frac{R_x}{nN}+\frac{1}{2},\\
&=\frac{\bar R_y-\bar R_x}{N}+\frac{1}{2}.
\end{align}

Next, the BM test considers a variance estimator for $\sqrt{N}\hat\theta$. If $\sqrt{N}$ is omitted, the estimated variance of $\hat\theta$ is written as follows:
\begin{align}
\mathrm{Var}(\hat\theta)&\simeq s^2=\frac{s_x^2}{nm^2}+\frac{s_y^2}{n^2m},\label{eq_bm_var}\\
s_x^2&=\frac{1}{n-1}\sum_{i=1}^n\left(R_i^x-r_i^x-\bar R_x+\frac{n+1}{2}\right)^2,\label{eq_sx}\\
s_y^2&=\frac{1}{m-1}\sum_{i=1}^m\left(R_i^y-r_i^y-\bar R_y+\frac{m+1}{2}\right)^2,\label{eq_sy}
\end{align}
where $s_x^2$ and $s_y^2$ are the sample variances of $R_i^x-r_i^x$ and $R_i^y-r_i^y$, respectively. Notice that $\bar r_x=(1/n)\sum_{i=1}^n r_i^x=(n+1)/2$ and $\bar r_y=(1/m)\sum_{i=1}^m r_i^y=(m+1)/2$. The meaning of (\ref{eq_bm_var}) is difficult to understand intuitively. I will explain later that the same expression can be derived using local linear approximation, providing a simpler interpretation. If $s_x^2=s_y^2=s^2=0$, which occurs when the two samples are completely separated as $\max(x_i)<\min(y_i)$ or $\max(y_i)<\min(x_i)$, these values may be replaced by small positive values, for example, $s^2_x=1/n$, $s^2_y=1/m$, and $s^2=2/(n^2m^2)$ \cite{schuurhuis2025}.

By using (\ref{eq_theta_hat2}) and (\ref{eq_bm_var}), the BM test statistic $W$ is written as follows:
\begin{align}
W&=\frac{\hat\theta-1/2}{s},\label{eq_w}\\
&=\frac{\bar R_y-\bar R_x}{N}\cdot\frac{nm}{\sqrt{ns_x^2+ms_y^2}},\\
&=\frac{nm(\bar R_y-\bar R_x)}{(n+m)\sqrt{ns_x^2+ms_y^2}}.
\end{align}

The test statistic $W$ asymptotically follows the standard normal distribution $\mathcal{N}(0,1)$ under the null hypothesis $\theta=1/2$. In practice, it is recommended to use a $t$-distribution whose degrees of freedom $\nu$ is calculated by the following equation \cite{brunner2000}.
\begin{equation}
\nu=\left.\left(\frac{s_x^2}{m}+\frac{s_y^2}{n}\right)^2 \middle/ \left(\frac{s_x^4}{(n-1)m^2}+\frac{s_y^4}{n^2(m-1)}\right)\right..\label{eq_df}
\end{equation}

The $p$-value for the two-tailed BM test can be calculated as follows:
\begin{equation}
p=2S_\nu(|W|),\label{eq_p}
\end{equation}
where $S_\nu$ is the survival function (1 minus the cumulative distribution function) of a $t$-distribution whose degrees of freedom is $\nu$.

There are several definitions of the 95~\% confidence interval associated with the BM test. The simplest one is as follows:
\begin{equation}
I_1=\hat\theta\ \pm\ Q_\nu (0.975) s,\label{eq_ci}
\end{equation}
where $Q_\nu$ is the quantile function (also called the percentile point function) of a $t$-distribution whose degrees of freedom is $\nu$. $Q_\nu\simeq 1.96$ for a large $\nu$. A drawback of this definition is that the upper bound may exceed 1 or the lower bound may become negative depending on the situation.

To overcome this, alternative definitions based on the logit transformation or the probit transformation have been proposed \cite{paulty2016}. The 95~\% confidence interval based on the logit transformation is written as follows:
\begin{equation}
I_2=g^{-1}\left(g(\hat\theta)\pm\frac{1.96s}{\hat\theta(1-\hat\theta)}\right),\quad g(x)=\log\frac{x}{1-x},\quad g^{-1}(x)=\frac{\exp(x)}{1+\exp(x)},\label{eq_ci2}
\end{equation}
where $g$ is the logit function. The 95~\% confidence interval based on the probit transformation is defined similarly.

\begin{figure}
\centering
\includegraphics[width=\linewidth]{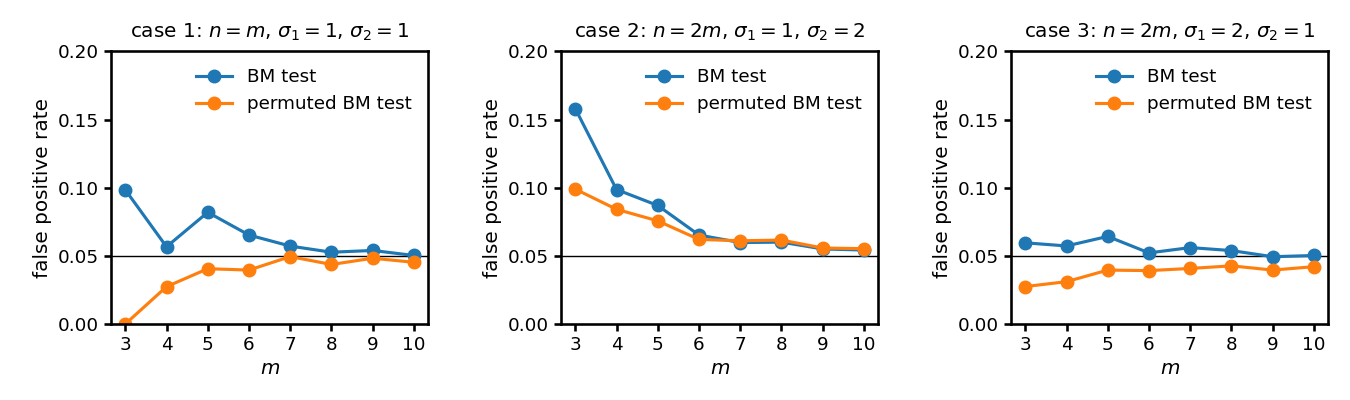}
\caption{Comparison of the Brunner-Munzel (BM) test and the permutation-based BM test for small sample sizes. Two independent samples $x_1,\ldots,x_n\sim\mathcal{N}(0,\sigma_1^2)$ and $y_1,\ldots,y_m\sim\mathcal{N}(0,\sigma_2^2)$ were considered. The significance level $\alpha$ was set to 0.05. The average false positive rate from 10{\thinspace}000 trials is shown. The maximum number of permutations was 10{\thinspace}000; this limit applied only when the total number of possible permutations exceeded that value.}
\label{fig7}
\end{figure}

If $\min(n,m)$ is small (e.g., less than 10), it was recommended to use a permutation test \cite{neubert2007}. Figure~\ref{fig7} shows the result of numerical simulations comparing the BM test \cite{brunner2000} and the permutation-based BM test \cite{neubert2007} for small sample sizes. Basically, the closer the false positive rate is to 0.05, the better the test is. However, upward (liberal) deviation from 0.05 is worse than downward (conservative) deviation in many practical situations because reporting statistical significance can have a greater impact on the target audience than reporting a non-significant result. The BM test tended to be liberal under the investigated conditions. The permutation-based BM test tended to be conservative in case 1 and case 3, and it was less liberal than the BM test in case 2. The irregular behavior observed in case 1 did not disappear when the number of trials was increased. These results suggest that the use of the permutation-based BM test should be considered if $\min(n,m)$ is small. 

\section{The \texorpdfstring{$C^2$}{C-square} test}

In this section, the $C^2$ test \cite{schuurhuis2025} is briefly introduced. The $C^2$ test is another nonparametric test for comparing two independent samples. One may call it the Schüürhuis-Konietschke-Brunner test or the SKB test to avoid confusion with the chi-square test. According to the original paper \cite{schuurhuis2025}, the $C^2$ test outperformed the BM test in broad conditions when $\min(n,m)\geq 15$. It was recommended not to use the $C^2$ test when $\min(n,m)<15$. Although the BM test is known to become inaccurate when the significance level $\alpha$ is $0.005$ \cite{noguchi2021}, the $C^2$ test remained accurate even when $\alpha=0.005$. Moreover, when $\alpha=0.05$, the $C^2$ test showed slightly better performance than the BM test. 

The $C^2$ test is based on an improved variance estimator for $\hat\theta$ that was recently proposed for $\min(n,m)\geq 2$ \cite{brunner2025b}:
\begin{align}
\tilde s^2&=\frac{(n-1)s_x^2+(m-1)s_y^2-nm\hat\theta(1-\hat\theta)+t_{xy}/4}{n(n-1)m(m-1)},\\
&=\frac{s_x^2}{nm(m-1)}+\frac{s_y^2}{n(n-1)m}-\frac{\hat\theta(1-\hat\theta)}{(n-1)(m-1)}+\frac{t_{xy}}{4n(n-1)m(m-1)},\label{eq_s2_tilde}
\end{align}
where $t_{xy}$ is the number of ties between $x_i$ and $y_j$:
\begin{equation}
t_{xy}=\sum_{i=1}^n\sum_{j=1}^m I(x_i=y_j),\quad I(x_i=y_j)=\left\{\begin{array}{ll}1&\text{if}\quad x_i=y_j,\\0&\text{if}\quad x_i\neq y_j.\end{array}\right.
\end{equation}

The first and second terms of (\ref{eq_s2_tilde}) resemble (\ref{eq_bm_var}). The third term corrects the overestimation in (\ref{eq_bm_var}). The fourth term is necessary only for discrete or ordinal data.

The new estimator $\tilde s^2$ is not only consistent but also unbiased. Its range is restricted as:
\begin{equation}
0\leq \tilde s^2\leq \frac{\hat\theta(1-\hat\theta)}{\min(n,m)-1}.
\end{equation}

The $C^2$ test statistic is defined as follows \cite{schuurhuis2025}:
\begin{equation}
C^2=\left(\frac{\hat\theta-1/2}{\tilde s}\right)^2 4\hat\theta(1-\hat\theta),\label{eq_c2}
\end{equation}
which approximately follows a $\chi^2$-distribution with one degree of freedom under the null hypothesis. If $\tilde s=0$, it is recommended to set $C^2=\min(n,m)$ \cite{schuurhuis2025}. Taking a square root of $C^2$ gives:
\begin{equation}
C=\frac{\hat\theta-1/2}{\tilde s}\cdot 2\sqrt{\hat\theta(1-\hat\theta)},\label{eq_c}
\end{equation}
which approximately follows the standard normal distribution $\mathcal{N}(0,1)$ under the null hypothesis. Although this formula is based on advanced analysis, the effect of the scale factor $2\sqrt{\hat\theta(1-\hat\theta)}$ can be intuitively understood. It is close to 1 when $\hat\theta\approx 1/2$, and it decreases as $\hat\theta$ deviates from $1/2$, reducing false positives. However, this does not reduce statistical power.

The $p$-value for the $C^2$ test can be calculated as follows:
\begin{equation}
p=2S(|C|),\label{eq_p_c2}
\end{equation}
where $S$ is the survival function of the standard normal distribution.

The 95~\% confidence interval associated with the $C^2$ test is calculated for $0<\hat\theta<1$ as follows \cite{schuurhuis2025}:
\begin{equation}
I_3=\frac{2\hat\theta+a\mp\sqrt{a^2+4a\hat\theta(1-\hat\theta)}}{2(1+a)},\quad a=1.96^2\frac{\tilde s^2}{\hat\theta(1-\hat\theta)}.\label{eq_ci3}
\end{equation}

It is guaranteed that $I_3\subseteq [0,1]$ holds. Moreover, $p< 0.05\Leftrightarrow 1/2\notin I_3$ is guaranteed, too.

Equations (\ref{eq_s2_tilde}) and (\ref{eq_ci3}) are complicated, and it may be difficult to understand their meanings intuitively. However, this may not necessarily prevent practitioners from using the $C^2$ test because it is likely that a simple interpretation is impossible. Indeed, the Welch's $t$-test is widely used even though the probability density function of the $t$-distribution is complicated.

The following sections focus on the interpretation of (\ref{eq_bm_var}), whose simple expression suggests that an intuitive explanation might be possible.

\section{Local linear approximation}

In this section, local linear approximation for estimating the variance of a statistic is explained. Let $X_\mathrm{all}=(X_1,\ldots,X_n)$ denote independent random variables, and let $x_\mathrm{all}=(x_1,\ldots,x_n)$ denote corresponding observations. The distributions followed by each $X_i$ are not necessarily identical. Let $T$ denote a statistic calculated from $X_\mathrm{all}$. Then, local linear approximation of $T$ for $X_\mathrm{all}\approx x_\mathrm{all}$ is written as follows:
\begin{equation}
T\simeq c+\sum_{i=1}^n g_i(X_i),
\end{equation}
where $c$ is a constant, and $g_1,\ldots,g_n$ are some functions. This means that the influence of each variable $X_i$ on the statistic $T$ is assumed to be additive, at least locally. If $g_1(X_1),\ldots,g_n(X_n)$ are approximately independent of each other for $X_\mathrm{all}\approx x_\mathrm{all}$, the variance of $T$ can be approximated for $X_\mathrm{all}\approx x_\mathrm{all}$ as follows:
\begin{align}
\mathrm{Var}(T)&\simeq \sum_{i=1}^n\mathrm{Var}(g_i(X_i)),\\
&\simeq\sum_{i=1}^n\mathrm{Var}(T\mid X_{-i}=x_{-i}),\label{eq_lla}
\end{align}
where $X_{-i}$ is $X_\mathrm{all}$ excluding $X_i$, and $x_{-i}$ is $x_\mathrm{all}$ excluding $x_i$. This corresponds to varying each variable one by one while keeping the other variables fixed.

For example, if $T=\bar X=(1/n)\sum_{i=1}^n X_i$, we get $\mathrm{Var}(T\mid X_{-i}=x_{-i})=\mathrm{Var}(X_i/n)=\mathrm{Var}(X_i)/n^2$ because only $X_i$ varies. Moreover, if all of $X_1,\ldots,X_n$ follow the same distribution with variance $\sigma^2$, we get $\mathrm{Var}(T)\simeq \sum_{i=1}^n\sigma^2/n^2=\sigma^2/n$, which is globally correct in this case.

This approximation is similar to the Efron-Stein inequality \cite{efron1981, steele1986}. However, it is an approximation rather than inequality, and it uses observations instead of expectations.

This method can be applied to the two-sample problem. Let $X_\mathrm{all}=(X_1,\ldots,X_n)$ and $Y_\mathrm{all}=(Y_1,\ldots,Y_m)$ denote independent random variables, and let $x_\mathrm{all}=(x_1,\ldots,x_n)$ and $y_\mathrm{all}=(y_1,\ldots,y_m)$ denote corresponding observations. Each $X_i$ follows the same distribution $F_x$, and each $Y_i$ follows the same distribution $F_y$. Let $T$ denote a statistic calculated from both $X_\mathrm{all}$ and $Y_\mathrm{all}$. Then, local linear approximation of $T$ for $X_\mathrm{all}\approx x_\mathrm{all}$ and $Y_\mathrm{all}\approx y_\mathrm{all}$ is written as follows:
\begin{equation}
T\simeq c+\sum_{i=1}^n g_i(X_i)+\sum_{i=1}^m h_i(Y_i),
\end{equation}
where $c$ is a constant, and $g_1,\ldots,g_n, h_1,\ldots,h_m$ are some functions. If $g_1(X_1),\ldots,g_n(X_n),h_1(Y_1),\ldots,h_m(Y_m)$ are approximately independent of each other for $X_\mathrm{all}\approx x_\mathrm{all}$ and $Y_\mathrm{all}\approx y_\mathrm{all}$, the variance of $T$ is approximated for $X_\mathrm{all}\approx x_\mathrm{all}$ and $Y_\mathrm{all}\approx y_\mathrm{all}$ as follows:
\begin{align}
\mathrm{Var}(T)&\simeq \sum_{i=1}^n\mathrm{Var}(g_i(X_i))+\sum_{i=1}^m\mathrm{Var}(h_i(Y_i)),\\
&\simeq\sum_{i=1}^n\mathrm{Var}(T\mid X_{-i}=x_{-i},Y_\mathrm{all}=y_\mathrm{all})+\sum_{i=1}^m\mathrm{Var}(T\mid X_\mathrm{all}=x_\mathrm{all},Y_{-i}=y_{-i}),\label{eq_lla2}
\end{align}
where $X_{-i}$ is $X_\mathrm{all}$ excluding $X_i$, $x_{-i}$ is $x_\mathrm{all}$ excluding $x_i$, $Y_{-i}$ is $Y_\mathrm{all}$ excluding $Y_i$, and $y_{-i}$ is $y_\mathrm{all}$ excluding $y_i$. 

For example, if $T=\bar Y-\bar X=(1/m)\sum_{i=1}^mY_i-(1/n)\sum_{i=1}^nX_i$, we get $\mathrm{Var}(T)\simeq \sum_{i=1}^n\mathrm{Var}(X_i/n)+\sum_{i=1}^m\mathrm{Var}(Y_i/m)\simeq \sigma_x^2/n+\sigma_y^2/m$, where $\sigma_x^2$ and $\sigma_y^2$ are the variances of $F_x$ and $F_y$, respectively. If $\sigma_x^2$ and $\sigma_y^2$ are replaced with respective sample variances, we get the variance estimator for $\bar Y-\bar X$ used in the Welch's $t$-test.

\section{Alternative interpretation of the Brunner-Munzel test}

In this section, a simple interpretation of the BM test based on the local linear approximation method is explained. Given observations $x_1,\ldots,x_n\sim F_x$ and $y_1,\ldots,y_m\sim F_y$, an $n\times m$ matrix $A=\{A_{ij}\}$ is defined as follows:
\begin{equation}
A_{ij}=\left\{\begin{array}{cc}1&\text{if}\quad x_i<y_j,\\1/2&\text{if}\quad x_i=y_j,\\0&\text{if}\quad x_i>y_j.\end{array}\right.
\end{equation}

Notice that the elements of $A$ are not independent of each other. For example, if $A_{11}=1$, $A_{21}=0$, and $A_{22}=1$ are given, that is, $x_1<y_1<x_2<y_2$ holds, $A_{12}$ cannot take a value of $0$ because $x_1<y_2$.

The sample stochastic superiority $\hat\theta$ corresponding to ({\ref{eq_theta_hat}}) is the mean of all elements of $A$, written as follows:
\begin{equation}
\hat\theta=\frac{1}{nm}\sum_{i=1}^n\sum_{j=1}^m A_{ij}.\label{eq_theta_hat3}
\end{equation}

For convenience, we also define the sum of each row and the sum of each column as follows:
\begin{align}
r_i&=\sum_{j=1}^m A_{ij},\quad (i=1,\ldots,n),\\
c_j&=\sum_{i=1}^n A_{ij},\quad (j=1,\ldots,m).
\end{align}

\begin{figure}
\centering
\includegraphics[width=0.8\linewidth]{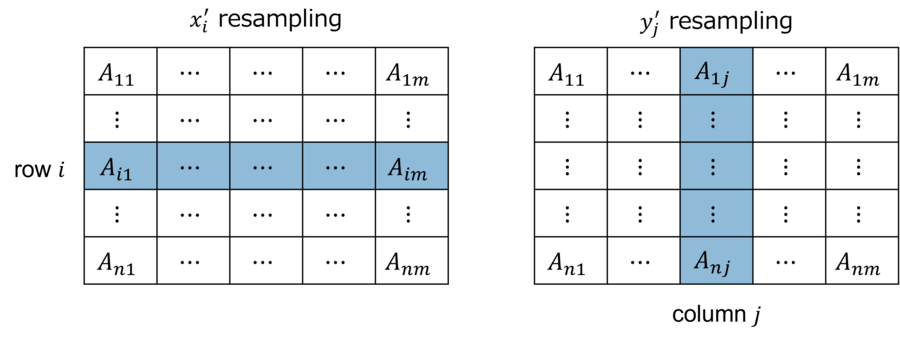}
\caption{Schematic of $x_i'$ resampling and $y_j'$ resampling. Shaded cells indicate cells updated through resampling of either $x_i'$ or $y_j'$.}
\label{fig8}
\end{figure}

Suppose that we can resample $x_i'\sim F_x$ for replacing $x_i$. This affects only the $i$-th row of $A$, as shown Fig.~\ref{fig8}. The sum of the $i$-th row $r_i$ changes to $r_i'$. Notice that the variance of $r_i'$---this will be obtained by repeating the resampling of $x_i'$ infinitely many times while keeping $y_1,\ldots,y_m$ fixed---can be approximated by the sample variance of $r_1,\ldots,r_n$. Similarly, if we can resample $y_j'\sim F_y$ for replacing $y_j$ while keeping $x_1,\ldots,x_n$ fixed, the sum of the $j$-th column $c_j$ changes to $c_j'$, whose variance can be approximated by the sample variance of $c_1,\ldots,c_m$.

\begin{figure}
\centering
\includegraphics[width=\linewidth]{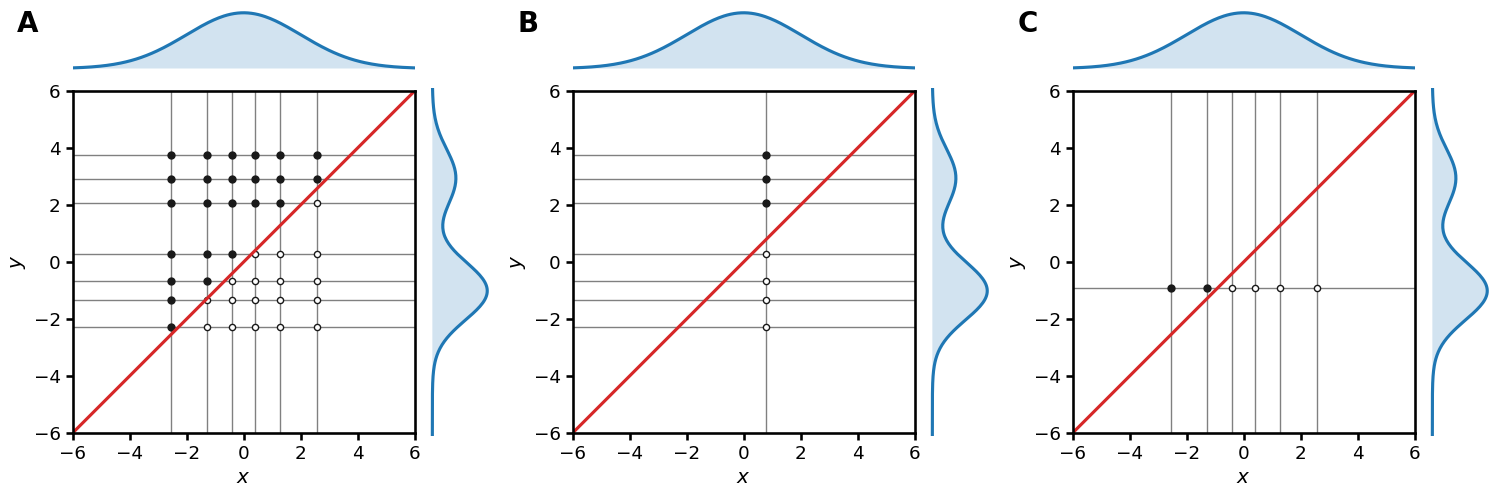}
\caption{Example of $x_i'$ resampling and $y_j'$ resampling. Red lines represent the line $y=x$. Vertical thin lines indicate the locations of $x$-observations. Horizontal thin lines indicate the locations of $y$-observations. Intersection points are shown in black when $y>x$, and in white otherwise. (A) All observations. (B) $x_i'$ resampling. (C) $y_j'$ resampling.}
\label{fig9}
\end{figure}

To clarify this point further, an example case is considered (Fig.~\ref{fig9}). There are $nm$ intersection points arising from all observations $x_1,\ldots,x_n$ and $y_1,\ldots,y_m$ (Fig.~\ref{fig9}A). Black points correspond to $A_{ij}=1$, and white points correspond to $A_{ij}=0$. When $x_i'$ is newly sampled from $F_x$ while keeping $y_1,\ldots,y_m$ fixed, $m$ intersection points arises (Fig.~\ref{fig9}B). Then, $r_i'$ can be understood as the number of black points on the thin vertical line. It varies depending on the location of $x_i'$ (in this example, it primarily takes a value of 3, 4, or 5). Its variance means the extent to which the number of black points fluctuates, and thus it can be approximated by the sample variance of $r_1,\ldots,r_n$. In this example, they are (7,5,4,3,3,2), as shown in Fig.~\ref{fig9}A. Similarly, when $y_j'$ is newly sampled from $F_y$ while keeping $x_1,\ldots,x_n$ fixed, $n$ intersection points arises (Fig.~\ref{fig9}C). Then, $c_j'$ can be understood as the number of black points on the thin horizontal line. Its variance is similarly approximated by the sample variance of $c_1,\ldots,c_m$. Iin this example, they are (1,1,2,3,5,6,6), as shown in Fig.~\ref{fig9}A.

Based on the above, the variance of $\hat\theta$ is estimated using local linear approximation ({\ref{eq_lla2}}) as follows:
\begin{align}
\mathrm{Var}(\hat\theta)&\simeq\sum_{i=1}^n\mathrm{Var}\left(\frac{r_i'}{nm}\right)+\sum_{j=1}^m\mathrm{Var}\left(\frac{c_j'}{nm}\right),\\
&=\frac{n\mathrm{Var}(r_1')+m\mathrm{Var}(c_1')}{n^2m^2},\\
&\simeq \frac{ns_x^2+ms_y^2}{n^2m^2}=s^2,\label{eq_bm_var2}\\
s_x^2&=\frac{1}{n-1}\sum_{i=1}^n\left(r_i-\bar r\right)^2,\label{eq_sx2}\\
s_y^2&=\frac{1}{m-1}\sum_{j=1}^m\left(c_j-\bar c\right)^2,\label{eq_sy2}
\end{align}
where $\bar r=(1/n)\sum_{i=1}^n r_i=m\hat\theta$ and $\bar c=(1/m)\sum_{j=1}^m c_j=n\hat\theta$. (\ref{eq_bm_var2}), (\ref{eq_sx2}), and (\ref{eq_sy2}) are equivalent to (\ref{eq_bm_var}), (\ref{eq_sx}), and (\ref{eq_sy}), respectively. The meaning of (\ref{eq_bm_var2}) is simple: the sum of variances obtained by assuming that each variable is virtually resampled one by one while keeping the other variables fixed.

The BM test statistic $W$ corresponding to (\ref{eq_w}) is defined as follows:
\begin{equation}
W=\frac{\hat\theta-1/2}{s},
\end{equation}
which can be calculated solely from $A$, without using in-group ranks and pooled ranks. Although ranks are computationally efficient, this matrix-based calculation can be performed in a very short time on modern computers anyway. A rank-free interpretation of the BM test similar to this has been proposed \cite{brunner2025}.

\section{Relation to jackknife resampling}

One might expect that the local linear approximation method described in Section~3 is related to jackknife (leave-one-out) resampling \cite{quenouille1956, tukey1958, miller1974}. The former is based on virtual resampling of each variable one by one, whereas the latter is based on removal of each variable one by one. In fact, Arvesen presented a generalized form of (\ref{eq_bm_var}) in the context of two-sample jackknifing \cite{arvesen1969}. However, its derivation was not explained in detail. Although (\ref{eq_bm_var}) is accurate unless $n$ and $m$ are small, which can be easily confirmed by numerical simulations, at least we can say that (\ref{eq_bm_var}) cannot be derived using jackknife resampling in a straight forward manner.

To explain this in more detail, I will briefly introduce the jackknife method and show the difficulty of applying it to estimate $\mathrm{Var}(\hat\theta)$.

\subsection{Jackknife resampling}

The jackknife method \cite{quenouille1956, tukey1958, miller1974} is a simple method for estimating the bias and variance of an estimator $\hat\theta$ for a parameter $\theta$. For example, if $\hat\theta$ is a sample mean, $\theta$ is a population mean. For simplicity, only the leave-one-out resampling is considered here. Let $\hat\theta$ denote an estimator for $\theta$  calculated from $n$ independent observations $x_1,\ldots,x_n\sim F_\theta$, where $F_\theta$ is a distribution depending on $\theta$. Let $\hat\theta_1,\ldots,\hat\theta_n$ denote \textit{jackknife replicates}, where each $\hat\theta_i$ is calculated from a subsample of $n-1$ observations omitting $x_i$. For example, $\hat\theta_1$ is calculated from $x_2,\ldots,x_n$, and $\hat\theta_2$ is calculated from $x_1,x_3,\ldots,x_n$.

Then, the following new variables are introduced, which are called \textit{pseudo-values}:
\begin{equation}
t_i=n\hat\theta-(n-1)\hat\theta_i,\quad (i=1,\ldots,n).
\end{equation}

For example, if $\hat\theta=\bar x=(1/n)\sum_{i=1}^n x_i$, we get $t_i=x_i$. In many cases, $t_1,\ldots,t_n$ are approximately independent of each other. Of course, there are cases where this does not apply (for example, the case of the median).

The jackknife estimate of $\theta$ is the average of them:
\begin{align}\hat\theta'&=\frac{1}{n}\sum_{i=1}^n t_i= n\hat\theta-(n-1)\bar\theta,\\
\bar\theta&=\frac{1}{n}\sum_{i=1}^n \hat\theta_i.\end{align}

This modification reduces the bias, if exists,  under several conditions \cite{quenouille1956}. For example, if $\hat\theta=(1/n)\sum_{i=1}^n(x_i-\bar x)^2$, it can be shown that $\hat\theta'=(1/(n-1))\sum_{i=1}^n(x_i-\bar x)^2$ \cite{quenouille1956}. The bias of $\hat\theta$ is approximated as follows:
\begin{equation}
E[\hat\theta]-\theta\simeq \hat\theta-\hat\theta'=(n-1)(\bar\theta-\hat\theta).
\end{equation}

If $t_1,\ldots,t_n$ are approximately independent of each other, the variance of $\hat\theta'$ can be approximated as follows:
\begin{align}
\mathrm{Var}(\hat\theta')&=\frac{1}{n^2}\mathrm{Var}\left(\sum_{i=1}^n t_i\right),\label{eq_jacvar_l1}\\
&\simeq \frac{1}{n^2}\sum_{i=1}^n \mathrm{Var}(t_i),\label{eq_jacvar_l2}\\
&= \frac{1}{n^2}\cdot n\mathrm{Var}(t_1),\label{eq_jacvar_l3}\\
&\simeq\frac{1}{n(n-1)}\sum_{i=1}^n (t_i-\bar t)^2,\label{eq_jacvar_l4}\\
&=\frac{n-1}{n}\sum_{i=1}^n(\hat\theta_i-\bar\theta)^2,\label{eq_jacvar_l5}
\end{align}
where $\bar t=(1/n)\sum_{i=1}^n t_i=\hat\theta'$. If $\hat\theta=\hat\theta'$, we get $V[\hat\theta]=V[\hat\theta']$. The independence assumption is used in the transformation from (\ref{eq_jacvar_l1}) to (\ref{eq_jacvar_l2}). It should be noted that the transformation from (\ref{eq_jacvar_l2}) to (\ref{eq_jacvar_l3}) is based on the assumption that $t_1,\ldots,t_n$ follow an identical distribution. Therefore, the final expression (\ref{eq_jacvar_l5}) is not directly applicable to two-sample cases, such as the BM test.

It can be seen that (\ref{eq_jacvar_l5}) is approximately $n$ times the sample variance of the jackknife replicates $\hat\theta_1,\ldots,\hat\theta_n$. Intuitively, when one observation is omitted, the value of $\hat\theta$ will change only slightly. However, if their influences are approximately independent, the sum of the changes resulting from removing each observation one by one is expected to be close to the change that would occur if all observations were resampled simultaneously. Such additivity of variances would be the key for the jackknife method to work correctly. For example, (\ref{eq_jacvar_l5}) becomes inaccurate when applied to the median due to nonlinear interaction between the jackknife replicates. Therefore, the fundamental idea behind jackknife resampling is similar to that of the local linear approximation method.

\subsection{Difficulty of estimating the variance of \texorpdfstring{$\hat\theta$}{q-hat} using jackknife resampling}

Unfortunately, it is difficult to estimate the variance of $\hat\theta$ defined in (\ref{eq_theta_hat3}) using jackknife resampling. Let $\hat\theta_1^x,\ldots,\hat\theta_n^x$ denote jackknife replicates, where each $\hat\theta_i^x$ is calculated from a subsample omitting $x_i$. Similarly, let $\hat\theta_1^y,\ldots,\hat\theta_m^y$ denote jackknife replicates, where each $\hat\theta_j^y$ is calculated from a subsample omitting $y_j$. Specifically, they can be written as follows:
\begin{align}
\hat\theta_i^x&=\frac{1}{(n-1)m}\sum_{\substack{k=1,\\k\neq i}}^n\sum_{j=1}^m A_{kj},\quad (i=1,\ldots,n),\\
\hat\theta_j^y&=\frac{1}{n(m-1)}\sum_{i=1}^n\sum_{\substack{k=1,\\k\neq j}}^m A_{ik},\quad (j=1,\ldots,m).
\end{align}

Then, pseudo-values are defined as follows \cite{arvesen1969}:
\begin{align}
t_i^x&=n\hat\theta-(n-1)\hat\theta_i^x=\frac{r_i}{m},\quad (i=1,\ldots,n),\\
t_j^y&=m\hat\theta-(m-1)\hat\theta_j^y=\frac{c_j}{n},\quad (j=1,\ldots,m).
\end{align}

The jackknife estimate of $\theta$ is defined as follow \cite{arvesen1969}:
\begin{equation}
\hat\theta'=\frac{1}{n+m}\left(\sum_{i=1}^n t_i^x+\sum_{j=1}^mt_j^y\right),
\end{equation}
but it turns out that $\hat\theta'=\hat\theta$. Therefore, $\hat\theta$ is already unbiased effectively; it does not contain the bias component that can be removed using the leave-one-out jackknife resampling.

Suppose that $t_1^x,\ldots,t_n^x, t_1^y,\ldots,t_m^y$ are approximately independent of each other. In other words, suppose that $r_1,\ldots,r_n,c_1,\ldots,c_m$ are approximately independent of each other. Then, the variance of $\hat\theta$ is calculated as follows:
\begin{align}
\mathrm{Var}(\hat\theta)&=\frac{1}{N^2}\mathrm{Var}\left(\sum_{i=1}^n t_i^x+\sum_{j=1}^m t_j^y\right),\\
&\simeq \frac{1}{N^2}\left(\sum_{i=1}^n \mathrm{Var}(t_i^x)+\sum_{j=1}^m \mathrm{Var}(t_j^y)\right),\\
&= \frac{1}{N^2}\left(n\mathrm{Var}(t_1^x)+m\mathrm{Var}(t_1^y)\right),\\
&= \frac{1}{N^2}\left(n\mathrm{Var}\left(\frac{r_1}{m}\right)+m\mathrm{Var}\left(\frac{c_1}{n}\right)\right),\\
&=\frac{1}{N^2}\left(\frac{n\mathrm{Var}(r_1)}{m^2}+\frac{m\mathrm{Var}(c_1)}{n^2}\right).\label{eq_bm_jacvar}
\end{align}

Notice that $\mathrm{Var}(r_1)$ is not conditional on $y_1,\ldots,y_m$, and thus it is much greater than $s_x^2$ that is calculated while keeping $y_1,\ldots,y_m$ fixed. Similarly, $\mathrm{Var}(c_1)$ is not conditional on $x_1,\ldots,x_n$ and is much greater than $s_y^2$. It is difficult to derive (\ref{eq_bm_var}) from (\ref{eq_bm_jacvar}) in a straight forward manner.

Taken together, it would be better to interpret the BM test from the viewpoint of local linear approximation rather than that of jackknife resampling.

\section{Conclusions}

In this paper, an alternative interpretation of the BM test statistic $W$ from the viewpoint of local linear approximation was introduced. Important points are as follows. (a) The sample stochastic superiority $\hat\theta$ is the average of $A_{ij}$, representing the balance between ones and zeros. (b) The variance estimator for $\hat\theta$ used in the BM test, denoted as $s^2$ in this paper, can be understood as the sum of variances obtained by assuming that each variable is virtually resampled one by one while keeping the other variables fixed. This approximation assumes local linearity. (c) The sum of each row $r_i$ and the sum of each column $c_j$ may help with interpretation. (d) Both $\hat\theta$ and $s^2$ can be calculated solely from $A$, without using in-group ranks and pooled ranks.

This simple interpretation will help practitioners decide to use the BM test without hesitation. Although the conditions under which local linear approximation holds are not clear, the accuracy of the BM test can be easily confirmed by numerical simulations. Given that the Welch's $t$-test is generally preferred to the Student's $t$-test as a parametric test for comparing two independent samples, the BM test and the $C^2$ test (the SKB test) will eventually supersede the $U$ test and become the standard non-parametric tests for comparing two independent sample \cite{karch2021}.

\section*{Acknowledgments}
The author wishes to express deep gratitude to Professor Edgar Brunner for giving numerous valuable comments and for pointing out the importance of the $C^2$ test.

\section*{Data availability}
The source code used for the numerical simulation in Fig.~6 is available in the GitHub repository \url{https://github.com/okumakito/bm2026} under the MIT license.

\section*{Funding}
This research was supported by JST Moonshot R\&D Grant Number JPMJMS2021.

\section*{Conflict of interest}
The author declares no competing interests.

\section*{Author contribution}
Makito Oku: Conceptualization, Methodology, Software, Formal analysis, Visualization, Writing---Original Draft, and Writing---Review \& Editing.

\section*{Artificial Intelligence tools}
No generative AI tools were used to write the manuscript. Only Google translation was used.

\bibliographystyle{ieeetr}
\bibliography{references}

\end{document}